\documentclass[9pt]{amsart}
\usepackage{amsmath}
\usepackage{amssymb}
\usepackage[all]{xy}

\numberwithin{equation}{section}

\newtheorem{defn}{Definition}[section]
\newtheorem{theorem}{Theorem}[section]
\newtheorem{assumption}[theorem]{Assumption}

\newtheorem{corollary}[theorem]{Corollary}

\newtheorem{prop}[theorem]{Proposition}
\newtheorem{remark}[theorem]{Remark}

\newtheorem{propdef}[theorem]{Proposition (Definition)}

\def\begineq{\begin{equation}}
\def\endeq{\end{equation}}

\def \n{\noindent}
\def \ggroup{\langle \mb g\rangle}

\def \mb{\mathbf}
\def \mc{\mathcal}

\def \bb{\mathbb}

\def \({\left(}
\def \){\right)}
\def \<{\langle}
\def \>{\rangle}
\def \xfrom{\xleftarrow}
\def \xto{\xrightarrow}

\def\qed{\hfill $\square$ \vspace{0.1in}}

\begin{document}
\title{A deRham model for Chen-Ruan cohomology ring of Abelian orbifolds}
\author{Bohui Chen}
\address{Department of Mathematics, Sichuan University, Chengdu, Sichuan 610064, P.R.CHINA}
\author{Shengda Hu}
\email{shengda@dms.umontreal.ca}
\address{
Centre de recherches math\'ematiques, Universit\'e de  Montr\'eal, CP 6128 succ Centre-Ville, Montr\'eal, QC H3C 3J7, Canada}

\abstract We present a deRham model for Chen-Ruan cohomology ring of abelian orbifolds. We introduce the notion of \emph{twist factors} so that formally the stringy cohomology ring can be defined without going through pseudo-holomorphic orbifold curves. Thus our model can be viewed as the classical description of Chen-Ruan cohomology for abelian orbifolds. The model simplifies computation of Chen-Ruan cohomology ring. Using our model, we give a version of wall crossing formula.
\endabstract
\maketitle

\section{Introduction}\label{intro}
In this paper we present a deRham model for Chen-Ruan cohomology ring of abelian orbifolds. We introduce the notion of \emph{twist factors} so that formally the stringy cohomology ring can be defined without going through pseudo-holomorphic orbifold curves. Thus our model can be viewed as the classical description of Chen-Ruan cohomology for abelian orbifolds. The model simplifies computation of Chen-Ruan cohomology ring and gives a version of wall crossing formula.

In their original papers \cite{CR1} and \cite{CR2}, the authors studied the Gromov-Witten theory of orbifolds. The theory in \cite{CR2} may be read as quantum cohomology ring theory of  orbifolds, while that in \cite{CR1}, as a special case of \cite{CR2}, serves as cohomology ring theory which is the now well-known Chen-Ruan cohomology ring of stringy orbifolds. We briefly review their construction for stringy abelian orbifolds in \S\ref{cr}.

The attempts of computing the Chen-Ruan cohomology ring structure is most successful for toric orbifolds and their hypersurfaces. The group structure of their Chen-Ruan cohomology is computed by M. Poddar \cite{P1}, \cite{P} and the Chen-Ruan ring structure for toric orbifolds is computed by Borisov et al,\cite{Bor}. In \cite{PP}, Parker et al computed the ring structure for the mirror quintic 3-fold. The difficulty of the computation in \cite{PP} stems from the fact that the Chen-Ruan cup product as defined in \cite{CR1} requires the computation of obstruction bundles over the moduli spaces of orbifold ghost curves.

In \S\ref{dr}, we propose a new formulation of Chen-Ruan cohomology for stringy abelian orbifolds. A deRham type theory is constructed with each cohomology class being represented by formal forms while the Chen-Ruan product is interpreted as ``wedge product" of formal forms. One may think of this as a classical level construction of Chen-Ruan theory. One advantage of the classical description is that it simplifies computations. To illustrate this point, in \S\ref{appl} we work out the computation of Chen-Ruan cohomology ring structure for the mirror quintic 3-fold and verifies the computations in \cite{PP}. Unfortunately, so far we have not found a similar way to deal with general orbifolds.

Let $G$ be a Lie group. The natural category for symplectic reduction with respect to Hamiltonian $G$ action is the category of symplectic orbifolds. As in the ordinary cohomology theory, it's natural to ask how the Chen-Ruan cohomology (ring) structure changes when crossing a wall. For example, the problem was posed in \cite{CR1}. Wall crossing have been studied by various authors for smooth cases. In \S\ref{wall}, we treat the problem for Chen-Ruan orbifold cohomology when $G$ is abelian. Our formulation leads to a natural extension of equivariant cohomology to $H^*_{G,CR}$ for torus action from which the surjectivity of the corresponding Kirwan map $\kappa: H^*_{G,CR} \to H_{CR}^*$ follows naturally. The main result is theorem \ref{prod:thm}, which is stated later. It reduces the change of the Chen-Ruan cohomology (ring) structure to computation at the fixed points in the wall. In \S\ref{appl} we apply the wall crossing formula to the simple case of weighted projective spaces, verifying the computation in \cite{J}.

Representing the cohomology class by forms in our new formulation (see \S\ref{dr}), we state the main theorem as following 

\n{\bf Theorem 4.3}
\emph{Let $G = S^1$ and $X$ be a Hamiltonian $S^1$-manifold with moment map $\mu : X \to \bb{R}$. Suppose $0 \in \bb{R}$ is a singular value and $F_{j \in J}$ be the fixed point components in $\mu^{-1}(0)$. Let $\tilde{\alpha}, \tilde{\beta}, \tilde{\gamma} \in H^*_{G,CR}(X)$ and $p, q \in \bb{R}$ be two regular values of $\mu$ such that $0 \in (p,q)$ is the only singular value. Denote $\alpha_p = \kappa_p(\tilde{\alpha})$ and so on, then we have}
$$\langle \alpha_q\cup \beta_q,\gamma_q\rangle
-\langle \alpha_p\cup \beta_p,\gamma_p\rangle
=
\sum_j
\int_{F_j} \frac{\tilde{i}_{(g_1)}(\tilde{\alpha}) \tilde{i}_{(g_2)}(\tilde{\beta})
\tilde{i}_{(g_3)}(\tilde{\gamma})}{e_G(N_{F_j})},$$
\emph{where $\cup$ is the Chen-Ruan cup product, $\langle\cdot, \cdot\rangle$ is the Poincar\'e pairing, $e_G(N_{F_j})$ is the equivariant Euler class of the normal bundle of $F_j$ in $X$ and $\tilde{i}_{(g)}(\cdot)$ is the equivariant twisted form defined by $\cdot$ and $g \in S^1$.}

\subsection*{Acknowledgements}
We'd like to thank Yongbin Ruan for posing the problem of wall crossing to 
us. The second author wants to thank Lev Borisov, Ernesto Lupercio and 
Yi Lin for helpful discussions. We also thank the authors of the excellent book \cite{ALR} for sharing with us the manuscript and the anonymous referee for insisting on using the more modern language of groupoids.

\section{Chen-Ruan cohomology theory for abelian orbifolds}
\label{cr}
In this section we review the theory of Chen-Ruan orbifold cohomology in the case where all the local isotropy groups are finite abelian groups. We refer to \cite{CR1}, \cite{CR2} and the excellent book \cite{ALR} for details and general setup.

\subsection{Abelian orbifolds}
\label{ab}
We recall briefly the language of groupoids and the definition of orbifolds in this language. 
A \emph{groupoid} $\mc G$ consists of the datum $(G_0, G_1; s, t, m, u, i)$ in the diagram:
$$\text{
\xymatrix{
{G_1 \left._s \times_t\right. G_1} \ar[r]^m & {G_1} \ar[r]^i & {G_1} \ar@<0.25ex>[r]^s \ar@<-0.25ex>[r]_t & {G_0} \ar[r]^u & {G_1}
}
},$$
where $G_0$ is the space of \emph{objects} and $G_1$ is the space of \emph{arrows}, with $s$ and $t$ being the \emph{source} and \emph{target} maps. The map $m$ defines \emph{composition} of two arrows while $i$ gives the \emph{inverse} arrow. The map $u$ is the \emph{unit} map, which is a two sided unit for the composition. The maps satisfies a set of natural axioms, such as $s(u(x)) = t(u(x)) = x$. We sometimes denote $i(g) = g^{-1}$ and $m(g, h) = gh$. The notion of \emph{morphism} between groupoids $\phi: \mc H \to \mc G$ consists of smooth maps $\phi_0 : H_0 \to G_0$ and $\phi_1 : H_1 \to G_1$ so that they are compatible with all the structure maps. Certain morphisms between groupoids are defined to be \emph{equivalences} and the \emph{Morita equivalence} between $\mc G$ and $\mc G'$ is defined by the existence of a groupoid $\mc H$ and the diagram $\mc G' \xfrom{\phi} \mc H \xto{\psi} \mc G$ where $\phi$ and $\psi$ are equivalences.

An \emph{orbifold groupoid} is defined to be a proper seperable \'etale Lie groupoid. It means that $G_0$ and $G_1$ are smooth Hausdorff manifolds and the structure maps are all smooth, with $s,t$ being local diffeomorphisms, so that $(s ,t): G_1 \to G_0 \times G_0$ is proper. It follows that $G_x = (s, t)^{-1}(x, x)$ for $x \in G_0$ is a finite group and is defined to be the \emph{isotropy} or \emph{local group} at $x$. The \emph{orbit space} $|\mc G|$ of $\mc G$ is defined to be the quotient space of $G_0$ under the equivalence relation $x\sim y$ iff they are connected by an arrow, i.e. $\exists g \in G_1$ so that $s(g) = x$ and $t(g) = y$. The simplest example for such a groupoid is the \emph{action groupoid} $G\ltimes M$ of a finite group $G$ acting on a manifold $M$, where $(G\ltimes M)_0 = M$ and $(G\ltimes M)_1 = G\times M$. The source and target maps are given by $(s, t) : (g, p) \mapsto (p, g\circ p)$ and the rest of structure maps are obvious. 

Suppose that $\phi: \mc G \to \mc H$ is an equivalence between orbifold groupoids, then induced map on the orbit spaces $|\phi| : |\mc G| \to |\mc H|$ is a homeomorphism. An \emph{orbifold structure} on a paracompact Hausdorff space $X$ is defined to be an orbifold groupoid $\mc G$ with a homeomorphism $f: |\mc G| \to X$ and $(\mc G, f)$ and $(\mc G', f')$ are \emph{equivalent} iff $\mc G$ and $\mc G'$ are Morita equivalent and the maps $f$ and $f'$ are compatible under the equivalence relation. 
Then an \emph{orbifold} $\mc X$ is defined to be a space $X$ with an equivalent class of orbifold structures. An orbifold structure $(\mc G, f)$ in such an equivalent class is a \emph{presentation} of the orbifold $\mc X$. We note that a presentation can be chosen such that over each point $x \in \mc X$, there is a component $\tilde U$ of $G_0$, so that the restricted groupoid is isomorphic to an action groupoid $G_{\tilde x}\ltimes \tilde U$ where $\tilde x \mapsto x$ under the quotient map. Such component $\tilde U$ is sometimes called an \emph{orbifold chart} around $x$. We thus see that, as abstract groups, $G_x$ is well defined for $x \in X$.
\begin{defn}
\label{ab:assump}
An orbifold $\mc X$ is an \emph{abelian orbifold} if the local groups $G_x$ for all $x \in X$ are abelian.
\end{defn}

\subsection{Twisted sectors}
\label{ab:twist}
We recall here the definition of twisted sectors in the language of groupoids. Let $\mc G$ be an orbifold groupoid, then a \emph{left $\mc G$-space} $M$ is a manifold with the \emph{anchor map} $\pi : M \to G_0$ and \emph{action map} $\mu : G_1 \left._s\times_\pi \right. M \to M$ satisfying the usual identies of an action:
$$\pi(\mu(g, p)) = t(g), \mu(u(x), p) = p \text{ and } \mu(g, \mu(h, p)) = \mu(m(g, h), p),$$
whenever the terms are well defined. Similar to the case of group actions, we may define the \emph{action groupoid} $\mc G \ltimes M$ for the $\mc G$-space $M$ with $(\mc G \ltimes M)_0 = M$ and $(\mc G \ltimes M)_1 = G_1 \left._s\times_\pi\right. M$. The source and target maps are given by $(s, t) : (g, p) \mapsto (p, \mu(g, p))$ as in the case of group actions. A special case is to let $M = G_0$, then the action groupoid is the groupoid $\mc G$.

In the following, we fix a groupoid presentation $\mc G$ of the orbifold $\mc X$ and will abuse notation and use $\mc G$ and $\mc X$ interchangeably. 
An \emph{orbifold morphism} is defined by a morphism between some groupoid presentations of the orbifolds. For any $x \in G_0$, it induces a map of isotropy groups $\phi_x : H_x \to G_x$. The morphism is called \emph{representable} if the map $\phi_x$ is injective for all $x \in G_0$. A morphism $\phi : \mc H \to \mc G$ of orbifold groupoids is an \emph{embedding} if $\phi_0$ is an immersion, $|\phi|$ is proper and satisfies a local condition which amounts to saying that the map $|\phi|$ can be locally lifted as a smooth map between the coverings. (cf. definition $2.3$ of the book \cite{SLR}). Then the pair $(\mc H, \phi)$ is a \emph{sub-groupoid} of $\mc G$ and correspondingly, the orbifold $\mc Y$ with underlying space $|\mc H|$ defined by $\mc H$ is a \emph{sub-orbifold} of $\mc X$. Sometimes, we abuse the notation and say that $\phi$, or $\mc H$ is a suborbifold. The \emph{intersection} of two suborbifold $\mc H$ and $\mc H'$ of $\mc G$ is defined to be the fibered product $\mc H \left._\phi\times_{\phi'}\right. \mc H'$ and will be denoted as usual $\mc H \cap \mc H'$.

The groupoid of \emph{twisted sectors $\wedge\mc G$} and \emph{$k$-multisectors $\mc G^k$} can be defined as action groupoid of certain left $\mc G$-space $\mc S_{\mc G}^k$ constructed naturally from $\mc G$:
$$\mc S_{\mc G}^k = \{(g_1, g_2, \ldots, g_k) \in G_1^k | s(g_1) = t(g_1) = s(g_2) = t(g_2) = \ldots = s(g_k) = t(g_k)\}.$$
The anchor map is $\pi_k : \mc S_{\mc G}^k \to G_0 : (g_1, g_2, \ldots, g_k) \mapsto x= s(g_1)$ and the action map for $s(h) = \pi_k (g_1, g_2, \ldots, g_k)$ is by conjugation $\mu_k(h, (g_1, g_2, \ldots, g_k)) = (hg_1h^{-1}, hg_2h^{-1}, \ldots, hg_kh^{-1})$. We note that when the orbifold is abelian, the action of $h \in G_x$ is trivial. 
In terms of the quotient orbifold, we have
$$\tilde X^k := |\mc G^k| = \{(x, (g_1, \ldots, g_k)_{G_x}) | x \in X, g_i \in G_x, i = 1, \ldots k\}.$$
Let $\mc S_o^k \subset \mc S_{\mc G}^k$ be the set of tuples so that $g_1g_2\ldots g_k = 1$, then it is a left $\mc G$-subspace and correspondingly defines a subgroupoid $\mc G_o^k$ of $\mc G^k$.

Associated to the orbifold groupoid $\mc G$, we have the \emph{skeletal groupoid} $\mc C$ with $C_i$ the discrete set of connected components of $G_i$, for $i = 0, 1$. The structure maps are induced from those of $\mc G$. Then $\mc C$ acts on the set $C(\mc S_{\mc G}^k)$ of connected components of $\mc S_{\mc G}^k$. Let $T^k = |\mc C \ltimes C(\mc S_{\mc G}^k)|$ and $(\mb g) = ((g_1, \ldots, g_k)_{G_x}) \in T^k$ the image of the component containing $(g_1, \ldots, g_k)$ with $s(g_1) = x$. Then $T^k$ parametrizes the connected components of $\tilde X^k$. We have the disjoint union of sub-groupoids
$$\mc G^k = \sqcup_{(\mb g) \in T^k} \mc G_{(\mb g)},$$
and correspondingly $\tilde {\mc X}^k = \sqcup_{(\mb g) \in T^k} \mc X_{(\mb g)}$. Analogously, let $T^k_o = |\mc C \ltimes C(\mc S_o^k)| \subset T^k$ and smiliarly define $\mc S_o^k$ and $\tilde{\mc X}^k_o$. Let $\mc S_{(\mb g)}$ be the preimage of $\mc X_{(\mb g)}$ under the natural quotient map, then it is a $\mc G$-subspace of $\mc S^k_{\mc G}$. The groupoid $\wedge \mc G$ is also called the \emph{inertia groupoid} of $\mc G$ and the corresponding orbifold $\wedge \mc X$ the \emph{inertia orbifold} or \emph{orbifold of twisted sectors} of $\mc X$. Correspondingly, $\mc X_{(\mb g)}$ is a \emph{$k$-multisector} or a \emph{twisted sector} when $k = 1$. The sector $\mc X_{(1)}$ associated to the unit is called the \emph{untwisted sector}. 
The multisectors are sub-orbifolds of $\mc X$, where the (union of) embedding(s) $\phi^k : \mc G^k \to \mc G$ is defined by $\phi^k_0 = \pi_k$ and $\phi^k_1(h, (g_1, \ldots, g_k)_{G_x}) = h$.
There are natural maps among the $k$-multisectors, which are induced by $\mc G$-equivariant maps among the $\mc S_{\mc G}^k$'s. The first class is the \emph{evaluation maps}, induced by
$$e_{i_1, \ldots, i_j} : \mc S_{\mc G}^k \to \mc S_{\mc G}^j : (g_1, g_2, \ldots, g_k) \mapsto (g_{i_1}, g_{i_2}, \ldots, g_{i_j}),$$
and the second class is the \emph{involutions}, induced by
$$I : \mc S_{\mc G}^k \to \mc S_{\mc G}^k : (g_1, \ldots, g_k) \mapsto (g_k^{-1}, \ldots, g_1^{-1}).$$
It's easy to see that the evaluation maps are (unions of) embeddings and the involutions are isomorphisms. In particular, the evaluation map $e = \pi_k: \mc S_{\mc G}^k \to G_0$ induces an embedding of $k$-multisectors as sub-orbifold (groupoid) of $\mc G$.

An \emph{orbifold bundle} $\mc E$ over $\mc G$ is by definition a $\mc G$-space so that $\pi: E \to G_0$ is a vector bundle and the action of $\mc G$ on $E$ is fiberwise linear, i.e. $g \in G_1$ induces linear isomorphism $g : E_{s(g)} \to E_{t(g)}$. The \emph{total space} $|\mc E|$ of $\mc E$ is given by the action groupoid $\mc G \ltimes E$ so that the projection morphism is defined by $\tilde \pi = (\pi, \pi_1)$ where $\pi_1 : G_1 \left._s\times_\pi\right. E \to G_1$ is the projection to the first factor. Then the map $|\tilde \pi| : |\mc E| \to |\mc G| = X$ gives the corresponding orbibundle. A \emph{section} $\sigma$ of bundle $\mc E$ is defined to be a $\mc G$-equivariant section of $E$, i.e. $s^*\sigma = t^*\sigma$ over $G_1$.
One example of orbifold bundle is the tangent bundle $T\mc G$, where the $\mc G$-space is $TG_0$ with the natural $\mc G$-action. The \emph{pull-back} $\phi^*\mc E$ of $\mc E$ by a morphism $\phi : \mc H \to \mc G$ is well-defined and is an orbifold bundle over $\mc H$.
Now suppose that $\phi: \mc H \to \mc G$ is an oriented suborbifold groupoid, then we define the \emph{normal bundle} of $\mc H$ in $\mc G$ as the quotient $N_{\mc H | \mc G} = \phi^*T\mc G / T\mc H$. 

\subsection{Degree shifting}
\label{ab:deg}
From now on, we assume that $\mc X$ is almost complex, that is, there is an almost complex structure on the tangent bundle $T G_0$, which is invariant under the $\mc G$-action. It follows that the $k$-multisectors have induced almost complex structures as well. 
The \emph{singular cohomology} $H^*(\mc X)$ of an orbifold $\mc X$ is defined to be the singular cohomology $H^*(X)$ of the underlying space $X = |\mc X|$. 
As ungraded group, the \emph{Chen-Ruan orbifold cohomology group} 
$H^*_{CR}(\mc X)$ of $\mc X$ is defined to be 
$$H^*_{CR}(\mc X) = H^*(\wedge \mc X) = \oplus_{(g) \in T_1}H^*(\mc X_{(g)}).$$
We now explain the grading. The degree of elements in $H^*(\mc X_{(g)})$ as elements in $H^*_{CR}(\mc X)$ is different from their degree in $H^*(\mc X_{(g)})$. The difference is the \emph{degree shifting number} $\iota_{(g)}$, which is defined below.

For $x \in \mc X_{(g)}$, let $g \in \mc S_{\mc G}$ be a preimage of $x$ and $\tilde x \in G_0$ be the image of $g$ under the evaluation map $e$. Then $g \in G_{\tilde x}$ and we have the decomposition into eigenspaces of $g$ action:
$$e^*T_{\tilde x} G_0 = T_g \mc S_{\mc G} \oplus N_g = T_g \mc S_{\mc G} \oplus \oplus_{j = 1}^m E_{j,g},$$
where $g$ action is trivial on the first summand and non-trivial on the rest. Choose a (complex) basis according to the above decomposition. Then the $g$ action can be represented by a diagonal matrix
\begin{equation}
\label{ab:diag}
diag (1, \ldots, 1, e^{2\pi i \theta_1},
 \ldots, e^{2 \pi i\theta_m}), \text{ where } 
\theta_j \in \mathbb{Q}\cap [0, 1) \text{ for all } j.
\end{equation}
The number $\iota(x, g) = \sum_{j}\theta_j$ doesn't depend on $x \in \mc X_{(g)}$ and is defined to be the \emph{degree shifting number} $\iota_{(g)}$ for the twisted sector $\mc X_{(g)}$.
In fact, the decompositions for all $x \in \mc X_{(g)}$ fit together and give a decomposition of tangent bundle with respect to $(g)$ action:
\begin{equation}
\label{ab:tandecom}
e^*T\mc G= T\mc X_{(g)} \oplus N_{(g)}
 = T\mc X_{(g)} \oplus {\oplus}_{j=1}^m E_j.
\end{equation}
We assume that each $E_j$ has rank $2$ where $2m$ is the codimension of $\mc X_{(g)}$. 
(In general, $E_j$ may not be a (complex) line bundle, in which case we may use standard splitting principle to proceed in the later arguments.) It's obvious that the bundle $N_{(g)}$ may be taken as the normal bundle $N_{\mc X_{(g)}|\mc X}$ of $\mc X_{(g)}$ in $\mc X$.

Under decomposition \eqref{ab:tandecom}, the matrix representing $g^{-1}$ can also be diagonalized:
$$diag (1, \ldots, 1, e^{2\pi i \theta'_1}, \ldots, e^{2 \pi \theta'_m}), 
\text{ where } \theta'_j \in \mathbb{Q}\cap [0, 1) \text{ for all } j.$$
Then it's easy to see that 
\begin{equation}
\label{complement}
\theta'_j + \theta_j = 
1 \Rightarrow \iota_{(g^{-1})} + \iota_{(g)} = m.
\end{equation}

Using the degree shifting number for the twisted sectors,
 we can write down $H^*_{CR}(\mc X)$ in graded pieces:
$$H^d_{CR}(\mc X) = \oplus_{(g)} H^d(\mc X_{(g)})[-2 \iota_{(g)}] = 
 \oplus_{(g)} H^{d-2\iota_{(g)}} (\mc X_{(g)}).$$
Note that, in general, the grading is rational instead of integral.

\subsection{Poincar\'e duality}
\label{ab:dual}
The Poincar\'e duality holds in \emph{de Rham} orbifold cohomology, which is isomorphic to the singular cohomology. The de Rham cohomology is defined as in the case of manifold, while the differential forms on an orbifold is defined to the the sections of the orbifold bundle $\wedge^* T^*\mc X$. Let $\alpha \in \Omega^*(\mc X)$ be a differential form, then the \emph{support} $\text{supp}(\alpha)$ \emph{of $\alpha$ in $\mc X$}  is the image in $X$ of its support in $G_0$. Let $U \subset G_0$ be an orbifold chart and let $\alpha$ be a differential form with compact support in $U$. Then the integration is defined by:
\begin{equation}
\label{orbint}
\int_U^{orb} \alpha = \frac{1}{|G_{U}|} \int_{U} \pi^*\alpha.
\end{equation}
Integration of general forms is then defined by partition of unity. For orbifolds admitting good covers, the pairing $\int^{orb}_{\mc X} \alpha_1 \wedge \alpha_2$ is a non-degenerate pairing between $H^*(\mc X)$ and $H^*_c(\mc X)$.

Poincar\'e duality holds in $H^*_{CR}(\mc X)$ with the involutions 
$$
I: \mc X_{(g)}\to \mc X_{(g^{-1})}.
$$
The pairing between $H^d_{CR}(\mc X)$ and $H^{2n-d}_{CR, c}(\mc X)$ is defined as the direct sum of the pairings on the twisted sectors $\mc X_{(g)}$ and $\mc X_{(g^{-1})}$:
$$\langle, \rangle^{(g)} : H^{d - 2\iota_{(g)}}(\mc X_{(g)}) \times H_c^{2n-d-2\iota_{(g^{-1})}} (\mc X_{(g^{-1})}) \to \mathbb{R} : \langle \alpha, \beta \rangle = \int^{orb}_{X_{(g)}} \alpha \wedge I^*(\beta).$$
Note that  $d - 2\iota_{(g)} + 2n-d-2\iota_{(g^{-1})} = 2(n-m)$ is the dimension of $\mc X_{(g)}$. The pairing $\langle, \rangle^{(g)}$ is simply the ordinary Poincar\'e duality on the abstract orbifold $\mc X_{(g)} \simeq \mc X_{(g^{-1})}$, which is non-degenerate.

\subsection{Obstruction bundles}
\label{ab:ob}
An important  ingredient
 in defining Chen-Ruan orbifold cup product
 is the \emph{obstruction bundles} on certain 
$3$-multisectors (also called \emph{triple twisted sectors}). Let $(\mathbf{g}) = ((g_1, g_2, g_3)) \in T^3_o$ and $E_{(\mathbf{g})} \to \mc X_{(\mathbf{g})}$ denote the obstruction bundle which we'll now describe.

Suppose that $\mathbf{r} = (r_1, r_2, r_3)$ records the orders of $g_i$ and let $(S^2, \mathbf{z}, \mathbf{r})$ be an orbifold $S^2$ with $3$ orbifold points $\mathbf{z} = (z_1, z_2, z_3)$. The local group at $z_i$ is the cyclic group of order $r_i$ for $i = 1, 2, 3$. Without loss of generality, we may assume that $\mathbf{z} = (0, 1, \infty)$ and drop $\mathbf{z}$ from the notation.

Fix an almost complex structure $J$ on the orbifold $\mc X$. We consider the space of representable pseudo-holomorphic orbifold morphisms $f: (S^2, \mathbf{r}) \to \mc X$, i.e. the local groups at $z_i$ are mapped injectively to the local groups of the image. In particular, we are interested in the maps where $[f] = 0 \in H_2(X)$, i.e. constant maps. The moduli space of such constant maps is given by $\tilde {\mc X}^3_o = \sqcup_{(\mb g) \in T^3_o} \mc X_{(\mb g)}$, which contains $\mc X_{(\mathbf{g})}$ as a connected component. The evaluation maps $e_i : \mc X_{(\mathbf{g})} \to \wedge \mc X$ for $i = 1, 2, 3$ play the same role as the usual evaluation maps on marked points.

Let $y = f(S^2) \in \mc X_{(\mathbf{g})}$, $e : \mc X_{(\mb g)} \to \mc X$ the evaluation map and consider the elliptic complex 
\begin{equation}
\label{orbidbar}
\bar{\partial}_y : \Omega^0((e\circ f)^*T\mc X) \to \Omega^{0,1}((e\circ f)^*T\mc X),
\end{equation}
which forms a family parameterized by $y \in \mc X_{(\mathbf{g})}$. The kernel of the family of elliptic complexes $\bar{\partial}_y$ is isomorphic to the bundle $T\mc X_{(\mathbf{g})}$ and the obstruction bundle $E_{(\mathbf{g})}$ is defined to be the cokernel.

More precisely,
let $\langle\mb g\rangle$ be the subgroup of $G_y$ generated by $\{g_i\}_{i = 1}^3$. It can be shown that $\<\mb g\>$, as abstract group, is independent of $y \in \mc X_{(\mb g)}$. Since $(\mb g) \in T^3_o$, there is a branched covering $\phi : \Sigma \to S^2$ from a smooth compact Riemann surface $\Sigma$ with covering group $\langle \mathbf{g} \rangle$ and branching loci over $(0, 1, \infty)$. The map $f: (S^2, \mathbf{r}) \to \mc X_{(\mathbf{g})}$ can then be lifted to the constant map $\tilde{f} : \Sigma \to \tilde{y} \in \mc S_{(\mb g)}$ where $\tilde{y} \mapsto y$ under the quotient.
Then the complex \eqref{orbidbar} lifts as the $\langle\mathbf{g}\rangle$-invariant part of the following complex
\begin{equation}
\label{liftdbar}
\bar{\partial}_{\tilde{y}} : \Omega^0((e\circ\tilde{f})^*TG_0) \to \Omega^{0,1}((e\circ\tilde{f})^*TG_0),
\end{equation}
where $e: \mc S_{(\mb g)} \to G_0$ is the evaluation map. Since $(e\circ\tilde{f})^*TG_0 = \Sigma \times T_{e(\tilde{y})}G_0$ is a trivial bundle, it follows that 
$$\mathrm{coker} \bar{\partial}_{\tilde{y}} = H^{0,1}(\Sigma) \otimes T_{e(\tilde{y})}G_0.$$
We see that the cokernel above fits together to give an orbifold bundle $H^{0,1}(\Sigma) \otimes e^*T\mc G$ over $\mc X_{(g)}$.
The group $\langle\mathbf{g}\rangle$ acts on the bundle with induced action on both factors. In particular, when the orbifold is abelian, the action of $G_y$ on the cokernel is given by the induced action on $T_{\tilde{y}}G_0$ and trivial action on $H^{0,1}(\Sigma)$.\footnote{In general, we should only take the action by $C(\mathbf{g})$, the elements in $G_y$ that commute with $\mathbf{g}$. Here $G_y$ is abelian, $C(\mathbf{g}) = G_y$.} Then the obstruction bundle is
$$E_{(\mb g)} = \left(H^{0,1}(\Sigma) \otimes e^*T\mc X\right)^{\langle\mathbf{g}\rangle}, \text{ where } e: \mc X_{(\mb g)} \to \mc X \text{ is the evaluation map}.$$
%

\subsection{Chen-Ruan orbifold cup product}
\label{ab:prod}

The definition of the Chen-Ruan orbifold cup product is the following. Let $(\mathbf{g}) \in T^3_o$ and $e_i : \mc X_{(\mathbf{g})} \to \mc X_{(g_i)}$ be the evaluation maps for $i = 1, 2, 3$. For $\alpha \in H^*(\mc X_{(g_1)})$, $\beta \in H^*(\mc X_{(g_2)})$ and $\gamma \in H^*_c(\mc X_{(g_3)})$, then $\alpha, \beta \in H^*_{CR}(\mc X)$ and $\gamma \in H^*_{CR, c}(\mc X)$, we define \emph{3-point function}
$$\langle \alpha, \beta, \gamma \rangle = \int^{orb}_{\mc X_{(\mathbf{g})}}e_1^*(\alpha) e_2^*(\beta) e_3^*(\gamma) e(E_{(\mathbf{g})}),$$
where $e(E_{(\mathbf{g})})$ is the Euler class of the obstruction bundle $E_{(\mathbf{g})}$ (computed by choosing a connection while the integral does not depend on the choice). The Chen-Ruan cup product $\alpha \cup \beta \in H^*(\mc X_{(g_3^{-1})}) \subset H^*_{CR}(\mc X)$ is then defined by Poincar\'e duality
\begin{equation}
\label{crproddef}
\langle \alpha \cup \beta, \gamma\rangle = \langle \alpha, \beta, \gamma \rangle \text{ for all } \gamma \in H^*_c(\mc X_{(g_3)}),
\end{equation}
It turns out that ``$\cup$'' defines an associative ring structure on $H^*_{CR}(\mc X)$.

\section{deRham model of $H^\ast_{CR}(\mc X)$} 
\label{dr}

It is well known from deRham theory that for a manifold $X$, the 
cohomology classes in $H^\ast(X)$ can be represented by
closed forms. 
In this section, we will present a similar model  
for $H^\ast_{CR}(\mc X)$ for abelian orbifold $\mc X$. Each cohomology class of (rational) degree $d$ will be represented by a formal $d$-form. A natural 
``wedge product" can be defined for these formal forms. We will show that 
this wedge product can be identified with Chen-Ruan orbifold cup product.
This somehow avoids the mysterious obstruction bundle.


\subsection{Twist factors}
\label{dr:fact}
To represent classes in 
$H^d(\mc X_{(g)})[-2\iota_{(g)}]$, besides 
a closed form on $\mc X_{(g)}$
we introduce an auxiliary  term
to account for the degree shifting.
The auxiliary term works as ``fractional Thom form''. We first recall the construction of Thom class in the category of orbifolds.
Let $\pi: \mc E \to \mc X$ be an oriented orbifold vector bundle, then a \emph{Thom form} $\Theta$ of $\mc E$ is defined by a $\mc G$-invariant Thom form of the vector bundle $E \to G_0$. As in the case of manifold, $\Theta$ is compactly supported on $\mc E$ near the $0$-section and we have for $\alpha \in \Omega^*(\mc E)$
\begin{equation}
\label{fact:notimp}
\int_{\mc E}^{orb} \alpha \wedge \Theta = \int_{\mc X}^{orb} i^*\alpha, \text{ where } i: \mc X \to \mc E \text{ is the } 0\text{-section}.
\end{equation}
The class represented by $\Theta$ is then the \emph{Thom class} of $\mc E$ and we denote it by $[\Theta]$.

Consider the splitting \eqref{ab:tandecom} and let $l_j$ be a Thom form of $E_j$ for $j = 1, \ldots, m$. Then the Thom class of $N_{(g)}$ is given by $\prod_{j=1}^m [l_j]$. We have
\begin{defn}
\label{dr:defact}The \emph{twist factor} $t(g)$ of $\mc X_{(g)}$ is defined by the \emph{formal}
 product
\begin{equation}
\label{fact:def}
[t(g)] = \prod_{j=1}^m [l_j]^{\theta_j},
\end{equation}
where $\theta_j$ is as given in \eqref{ab:diag}.
\end{defn}
Because of the similarity of \eqref{fact:def} with the Thom class, $[t(g)]$ may be regarded as \emph{fractional Poincar\'e dual} of $\mc X_{(g)}$ in $\mc X$.
Formally, $t(g)$ is a form of degree $2 \iota_{(g)}$ supported in a neighbourhood of $\mc X_{(g)}$ in $\mc X$. This formal degree
makes up the difference between the degrees of the classes in $H^*_{CR}(\mc X)$
and $H^*(\mc X_{(g)})$. We can then write the identification of 
$H^*(\mc X_{(g)})$ as a summand in  $H^*_{CR}(\mc X)$ formally as a 
 {\em Thom isomorphism}.
More precisely, suppose $U$ be a neighbourhood of $\mc X_{(g)}$ in its normal bundle and identify it with a neighbourhood of $\mc X_{(g)}$ in $\mc X$,
with the projection map $\pi: U\to \mc X_{(g)}$
 and by representing cohomology classes by forms
, we formally write
\begin{equation}
\label{dr:thomiso}
i_{(g)}: 
H^*(\mc X_{(g)}) \to H^{*+ 2\iota_{(g)}}_{CR}(\mc X) 
:[\alpha]\to [\pi^*(\alpha)t(g)].
\end{equation}
We shall call $i_{(g)}(\alpha)$ a \emph{twisted form} and note that it is supported in the neighbourhood $U$ of $\mc X_{(g)}$. We'll drop $\pi^*$ and simply write $i_{(g)}(\alpha) = \alpha t(g)$ since there should be no confusion.
\begin{remark}
\label{fact:extra}\rm{
In what follows, we'll carry out the (formal) computations using notations $[l_j]$ instead of $l_j$ to emphasis that the results in cohomology do not depend on the particular choice of the forms $l_j$. 
}\end{remark}

\subsection{Poincar\'e duality}
\label{dr:dual}
We will discuss the wedge product of two twisted forms later. As a warm up, we explain how
the Poincar\'e duality (\S\ref{ab:dual}) follows from this formulation.
We follow the convention that integration of the form $\int_{\mc X}^{orb} \alpha \wedge \prod_{j=1}^k t(g_j)$ vanishes unless the product $\prod_{j=1}^k [t(g_j)]$ gives the Thom class of some suborbifold $\mc Y$ of $\mc X$, in which case \eqref{fact:notimp} applies. 
Let $a=i_{(g)}(\alpha)$ and $b=i_{(g^{-1})} (\beta)$ and define pairing as
$$
\langle a,b\rangle
= \int ^{orb}_{\mc X} a\wedge b
=\int^{orb}_{\mc X} \alpha\wedge \beta\wedge t(g)\wedge t(g^{-1})
=\int^{orb}_{\mc X_{(g)}} \alpha\wedge \beta.
$$
For the last equality, we use the fact that $t(g)\wedge t(g^{-1})$
is the Thom form of $\mc X_{(g)}$ in $\mc X$ (cf. \eqref{complement}, \eqref{fact:notimp}).
This matches $\langle \alpha,\beta\rangle^{(g)}$ as defined in \S\ref{ab:dual}.

\subsection{Wedge product}
\label{dr:prod}
Now we assume that the orbifold $\mc X$ is abelian. Then the normal bundle $N_{\mc X_{(\mb g)}| \mc X}$ of any $k$-multisector $\mc X_{(\mb g)}$ can be decomposed into direct sum of (complex) line bundles with respect to the $\<\mb g\>$-action (upto splitting principle):
\begin{equation}\label{dr:multidecom}
N_{\mc X_{(\mb g)}| \mc X} = \oplus_{j} E_j.
\end{equation}
Let $a_i\in H^{d_i+2\iota_{(g_i)}}_{CR}(\mc X),i=1,2,$ be two twisted forms. The wedge product $a_3 = a_1\wedge a_2$
can be 
defined formally in the obvious way.
 We explain that $a_3$ is also a twisted form in a very natural way. 
\begin{propdef}
\label{prod:ring}
Suppose 
$a_i= i_{(g_i)}\alpha_i= \alpha_i t{(g_i)}, \text{ for } i = 1, 2$ 
and define $a_3 := a_1 \wedge a_2 = i_{(g_1)}(\alpha_1) \wedge i_{(g_2)}(\alpha_2)$.
Then $\exists \alpha_{3,j} \in H^*(\mc X_{(g_{3,j})})$ where $(g_1, g_2, g_{3,j}^{-1}) \in T^3_o$, so that $a_3 = \sum_{j}\alpha_{3,j} t(g_{3,j})\in H^*_{CR}(\mc X)$.
\end{propdef}
\n
{\bf Proof.}
The formula defining $a_3$ is interpreted as following. The form $t{(g_i)}$, and thus $a_i$, is supported near $\mc X_{(g_i)}$ for $i = 1,2$. It follows that $a_3$ is supported near $\mc Z= \mc X_{(g_1)}\cap \mc X_{(g_2)}$. We see that $\mc Z$ is a union of $2$-multisectors $\mc Z_j : = \mc X_{(h_{1,j}, h_{2,j})}$ so that $h_{i,j}$ is in $(g_i)$ for all $j$. Let $g_{3,j} = h_{1, j}h_{2,j}$, then $((h_{1,j}, h_{2,j}, g_{3,j}^{-1})) \in T^3_o$ and we have $\mc Z_j \xto{i_{3,j}} \mc X_{(g_{3, j})}$ is naturally an embedding.
We define
\begin{equation}
\label{prod:def}
a_{3,j} = \(a_1 \wedge a_2\)_j := \(i_{1,j}^*(\alpha_1) \wedge i_{2,j}^*(\alpha_2)\) \wedge t{(g_1)}\wedge t{(g_2)},
\end{equation}
where $i_{\bullet,j}$ is the inclusion $\mc Z_j \to \mc X_{(g_{\bullet})}$. Then it's obvious that $a_3 = \sum_j a_{3,j}$. Thus we only need to prove the proposition for the case where $\mc Z$ has only one component. The main issue is to deal with $t{(g_1)}\wedge t{(g_2)}$.

Assume that $\mc Z$ has only one component. It is clear that the normal bundle $N_{\mc Z|\mc X}$ of $\mc Z$ in $\mc X$ has the following splitting
$$
N_Z= N_1\oplus N_2\oplus N_{3}\oplus N',
$$
where $N_i = N_{\mc Z | \mc X_{(g_i)}}$ are the normal bundles of $\mc Z$ in $\mc X_{(g_i)}$ and $N'$ is defined by the equation.
$N_i,i=1,2,3$ and $N'$ are further decomposed into line eigenbundles
$$
N_i= \oplus_{j=1}^{k_i} L_{ij}; 
N'= \oplus_{j=1}^k L'_j.
$$
The splitting of normal bundle $N_{\mc X (g_i) | \mc X}$ restricting to $\mc Z$ is 
compatible with this splitting. For instance,
$$
N_{\mc X (g_1) | \mc X}|_Z=
N_2\oplus N_3\oplus N'=\left(\oplus_{j=1}^{k_2} L_{2j}\right)\bigoplus\left(\oplus_{j=1}^{k_3} L_{3j}\right)\bigoplus\left(\oplus_{j=1}^k L'_j\right)
$$
and correspondingly, near $\mc Z$ 
$$
t{(g_1)}= t_2(g_1)t_3(g_1)t'(g_1).
$$
The terms in the right hand side above are defined in the obvious way.
Similarly,
$$
t{(g_2)}= t_1(g_2)t_3(g_2)t'(g_2) \text{ and } t{(g_3)}= t_1(g_3)t_2(g_3)t'(g_3).
$$
We look at the \emph{formal} expression
\begin{equation}
\label{prod:factor}
\frac{t(g_1)\wedge t(g_2)}{t(g_3)}
=\frac{t_2(g_1)t_1(g_2)}{t_1(g_3)t_2(g_3)}\{t_3(g_1)t_3(g_2)\}\frac{t'(g_1)t'(g_2)}{t'(g_3)}.
\end{equation}
\begin{enumerate}
\item Note that $(g_1, g_2, g_3^{-1}) \in T^3_o$, then it's easy to check that the first fraction simplifies to $1$ when restricted to $\mc Z$.
\item
In light of the equations \eqref{fact:notimp} and \eqref{fact:def}, 
the term $\{t_3(g_1)t_3(g_2)\}$
 is a Thom form $\tau_3$ of $N_3$,
representing the Poincar\'e dual of $\mc Z$ in $\mc X_{(g_3)}$.
\item
To see what happens to $\frac{t'(g_1)t'(g_2)}{t'(g_3)}$, let us  look at each $L'_j$. Suppose
$g_i$ acts on $L'_j$ as multiplication of $e^{2\pi i\theta_{ij}}$ and the Thom form
of $L'_j$ is $[l'_j]$, then
the exponent of $[l'_j]$ in $t{(g_i)}$ is $\theta_{ij}$. Now $g_3 = g_1g_2$ implies that 
$\theta_{1j}+\theta_{2j}$ is either $\theta_{3j}$ or $\theta_{3j}+1$.
Hence,
$$
\frac{[l'_j]^{\theta_{1j}}[l'_j]^{\theta_{2j}}}
{[l'_j]^{\theta_{3j}}}
=\left\{
\begin{array}{ll}
1, & \mbox{if } \theta_{1j}+\theta_{2j}=\theta_{3j},\\
{[}l'_j{]}, & \mbox{if } \theta_{1j}+\theta_{2j}=\theta_{3j}+1.
\end{array}
\right.
$$
The right hand side of the above becomes ordinary forms when restricted to $\mc Z$. Set
$$\Theta_{(g_1,g_2)}=\left.\frac{t'(g_1)t'(g_2)}{t'(g_3)}\right|_{\mc Z},$$
then $\Theta_{(g_1,g_2)}$ is an ordinary form. In fact, it's quite clear that $[\Theta_{(g_1, g_2)}] = e\left(E'_{(g_1,g_2)}\right) \in H^*(\mc Z)$ where 
\begin{equation}
\label{prod:obfact}
E'_{(g_1,g_2)}
=\bigoplus_{\theta_{1j}+\theta_{2j}=\theta_{3j}+1} L'_j.
\end{equation}
\end{enumerate}
It follows that \eqref{prod:factor} gives an honest form on $\mc X_{(g_3)}$ and $a_3= i_{(g_3)}(\alpha_3)$ where
$$[\alpha_3]
=[\(i_1^*(\alpha_1)\wedge i_2^*(\alpha_2)\wedge \Theta_{(g_1,g_2)}\)\wedge 
\tau_3]
 \in H^*(\mc X_{(g_3)}),
$$
is given by $[i_1^*(\alpha_1)\wedge i_2^*(\alpha_2)\wedge \Theta_{(g_1,g_2)}] \in H^*(\mc Z)$ via Thom homomorphism. \qed
\begin{corollary}
\label{dr:ringstr}
The product $\wedge$ defines an associative ring structure on $H^*_{CR}(\mc X)$.
\end{corollary}
\n
{\bf Proof.} Assuming all intersections in the following has only one component. The general case is dealt with similarly as in the proposition. The equation we need to establish is:
\begin{equation*}
\begin{split}
(a_1 \wedge a_2) \wedge a_3 & = i_{(g_5)}\left\{\(i_4^*\left\{(i_1^*(\alpha_1) \wedge i_2^*(\alpha_2) \wedge \Theta_{(g_1, g_2)}) \wedge \tau_4\right\} \wedge i_3^*(\alpha_3) \wedge\Theta_{(g_4, g_3)} \) \wedge \tau_5\right\} \\
= a_1 \wedge(a_2\wedge a_3) & = i_{(g_5)}\left\{\(i_1^*(\alpha_1) \wedge i_{4'}^*\left\{(i_2^*(\alpha_2) \wedge i_3^*(\alpha_3) \wedge \Theta_{(g_2, g_3)}) \wedge \tau_{4'}\right\} \wedge \Theta_{(g_1, g_{4'})}\)  \wedge \tau_5 \right\},
\end{split}
\end{equation*}
where on the left hand side $\mc X_{(g_1)} \cap \mc X_{(g_2)} = \mc Z_4 \subset \mc X_{(g_4)}$ with $(g_1, g_2, g_4^{-1}) \in T^3_o$ and $\mc X_{(g_4)} \cap \mc X_{(g_3)} = \mc Z_5 \subset \mc X_{(g_5)}$ with $(g_4, g_3, g_5^{-1}) \in T^3_o$ and on the right hand side $\mc X_{(g_2)} \cap \mc X_{(g_3)} = \mc Z_{4'} \subset \mc X_{(g_{4'})}$ with $(g_1, g_2, g_{4'}^{-1}) \in T^3_o$ and $\mc X_{(g_{4'})} \cap \mc X_{(g_3)} = \mc Z_5 \subset \mc X_{(g_5)}$ with $(g_{4'}, g_3, g_5^{-1}) \in T^3_o$. The notation "$\subset$" denotes embedding given by the composition of arrows. The rest of the notations are as in the proposition. Let $\mc Z = \mc X_{(g_1)} \cap \mc X_{(g_2)} \cap \mc X_{(g_3)}$ then both sides of the equation is supported in a neighbourhood of $\mc Z$. We rewrite the left hand side:
$$LHS = i_{(g_5)}\left\{\(\(i_1^*(\alpha_1) \wedge i_2^*(\alpha_2) \wedge i_3^*(\alpha_3)\) \wedge \(i_4^*(\Theta_{(g_1, g_2)}) \wedge \Theta_{(g_4, g_3)}\)\) \wedge \(i_4^*(\tau_4)\wedge \tau_5\)\right\} $$
where $i_{\bullet}$ is the inclusion of $\mc Z$ in $\mc X_{(g_{\bullet})}$ for $\bullet = 1, 2, 3$ and $i_4$ is the inclusion of $\mc Z_5$ in $\mc X_{(g_4)}$. Then $\(i_4^*(\tau_4)\wedge \tau_5\) = \tau_{\mc Z}(\mc X_{(g_5)})$. It follows from \eqref{prod:obfact} that $\(i_4^*(\Theta_{(g_1, g_2)}) \wedge \Theta_{(g_4, g_3)}\)$ represents the Euler class of the following bundle over $\mc Z$:
$$E'_{(g_1, g_2, g_3)} = \bigoplus_{\theta_{1j}+\theta_{2j} + \theta_{3j}=\theta_{5j}+2} E_j,$$
where $E_j$'s are the complex line bundles appearing in the decomposition \eqref{dr:multidecom} for $\mc X_{(g_1, g_2, g_3)}$. Thus we have:
$$LHS = i_{(g_5)}\left\{\(\(i_1^*(\alpha_1) \wedge i_2^*(\alpha_2) \wedge i_3^*(\alpha_3)\) \wedge e\(E'_{(g_1, g_2, g_3)}\)\) \wedge \tau_Z(X_{(g_5)})\right\}$$
The equation can then be shown by similar rewriting of the right hand side.
\qed

\subsection{Obstruction bundle and obstruction form}\label{dr:ob}
The natural map $T^2 \to T^3_o : ((g, h)) \mapsto ((g, h, (gh)^{-1}))$ is an isomorphism and we have isomorphism $\tilde{\mc X}^2 \to \tilde{\mc X}^3_o$ correspondingly. We consider a component $\mc Z_j = \mc X_{((h_{1, j}, h_{2,j}))}$ as in the proof of proposition \ref{prod:ring} and use $\mc Z_j$ to denote also the $3$-multisector $\mc X_{((h_{1, j}, h_{2,j}, g_{3,j}^{-1}))}$, in the notation of the previous section. We show that
\begin{prop}
$\Theta_{(h_{1,j},h_{2,j})}= e\left(E_{((h_{1, j}, h_{2,j}, g_{3,j}^{-1}))}\right)$ on $\mc Z_j$.
\end{prop}
\n
{\bf Proof.}
As we are only considering one component, we let $g_1 = h_{1, j}$, $g_2 = h_{2, j}$, $g_3 = g_{3,j}$ and $\mc Z = \mc Z_j$. 
It then suffices to show that $E'_{(g_1,g_2)}\cong E_{(\mb g)}$ on $\mc Z$. We'll use the notations in \S\ref{ab:ob}.

Let $e: \mc Z \to \mc X$ be the evaluation map. With decomposition \eqref{dr:multidecom} for $\mc Z$
and the almost complex structure on $\mc X$, 
the matrices representing the action of elements in $\langle \mathbf{g} \rangle$ can all be diagonalized. In particular 
we have
$$g_i = diag(1, \ldots, 1, e^{2\pi i \theta_{i1}}, \ldots, e^{2\pi i\theta_{im}}), \text{ where } \theta_{ij} \in \mathbb{Q}\cap [0, 1), \text{ for } i = 1, 2, 3.$$
The fiber of $E_{(\mathbf{g})}$ at $y$ is then
\begin{equation}
\label{prod:ob}
\begin{split}
E_{(\mathbf{g}),y} = & (H^{0, 1}(\Sigma) \otimes T_{e(\tilde{y})}G_0)^{\langle \mathbf{g}\rangle} \\
=& 
(H^{0, 1}(\Sigma) \otimes T_{\tilde{y}}\mc S_{(\mb g)})^{\langle \mathbf{g}\rangle}\oplus
\oplus_{j=1}^m
(H^{0, 1}(\Sigma) \otimes E_j|_y)^{\langle \mathbf{g}\rangle} \\
= & H^1\(S^2, \(\phi_*(T_{\tilde{y}}\mc S_{(\mb g)})\)^{\langle \mathbf{g}\rangle}\) \oplus \oplus_{j=1}^m H^1\(S^2, \(\phi_*(E_j|_y)\)^{\langle \mathbf{g}\rangle}\)
\end{split}
\end{equation}
where $\phi : \Sigma \to S^2$ is the branched covering, $\phi_*$ is the push-forward of the constant sheaves. Let $V$ be $\ggroup$-vector space of (complex) rank $v$ and $m_{i,j} \in \bb{Z}\cap [0, r_i)$ be the weights of action of $g_i$ on $V$. Applying the index formula (proposition 4.2.2 in \cite{CR1})
to $\(\phi_*(V)\)^{\langle \mathbf{g}\rangle}$ we have
$$
\chi=
v-\sum_{i=1}^3\sum_{j = 1}^v\frac{m_{i,j}}{r_i},
$$
 Here we use the fact that $c_1(\phi_*(V))=0$ for constant sheaf $V$. If $\ggroup$ action is trivial on $V$ then $\chi = v$. For $V=E_j|_y$, we see that $v = 1$ and $\frac{m_{i,1}}{r_i}$ is just 
$\theta_{ij}$.

With the above preparations, we have the following
\begin{enumerate}
\item
$(H^{0, 1}(\Sigma) \otimes T_{\tilde{y}}\mc S_{(\mb g)})^{\langle \mathbf{g}\rangle}=\{0\}$ and \\
\item
$(H^{0, 1}(\Sigma) \otimes E_j|_y)^{\langle \mathbf{g}\rangle}$
is nontrivial ($\Rightarrow$ rank $1$) $\iff \sum_{i=1}^3 \theta_{ij}=2.$ (Note that this sum is either 1 or 2.)
Moreover, 
it is clear that
$$
(H^{0, 1}(\Sigma) \otimes E_j|_y)^{\langle \mathbf{g}\rangle}
\cong E_j|_y.
$$
\end{enumerate}
It follows that
$$E_{(\mathbf{g})} = \bigoplus_{\sum_{i = 1}^3\theta_{ij} = 2} E_j,$$
which obviously matches with \eqref{prod:obfact} (since the $\theta_{3j}$ here is $1-\theta_{3j}$ in \eqref{prod:obfact}). \qed

\subsection{Ring isomorphism}
\label{dr:iso}
So far, on $H^\ast_{CR}(\mc X)$ we have two different product structures: Chen-Ruan product 
``$\cup$'' and wedge product ``$\wedge$''.
We have 
\begin{theorem}
\label{isothm}
$(H^\ast_{CR}(\mc X), \cup)\cong (H^\ast_{CR}(\mc X),\wedge)$ as rings.
\end{theorem}
\noindent
{\bf Proof.}
Let $\alpha,\beta$ and $\gamma$ be as in \S \ref{ab:prod}.
We show that 
$$
\langle\alpha\cup \beta,\gamma\rangle= 
\int_{\mc X}^{orb} i_{(g_1)}(\alpha)\wedge i_{(g_2)}(\beta)
\wedge i_{(g_3)}(\gamma).
$$
The right hand side is
\begin{equation*}
\begin{split}
\int_{\mc X}^{orb} i_{(g_1)}(\alpha) i_{(g_2)}(\beta) i_{(g_3)}(\gamma) = & \int_{\mc X}^{orb} \pi_1^*(\alpha) \pi_2^*(\beta) \pi_3^*(\gamma) \prod_{i = 1}^3 t{(g_i)} \\
= & \int_{\mc X}^{orb} \pi_1^*(\alpha) \pi_2^*(\beta) \pi_3^*(\gamma) \prod_{j=1}^m [l_j]^{\sum_{i = 1}^3\theta_{ij}} \\
= & \int_{\mc X}^{orb}
 \pi_1^*(\alpha) \pi_2^*(\beta) \pi_3^*(\gamma) \Omega(X_{(\mb g)}) \prod_{j = 1}^m [l_j]^{\sum_{i = 1}^3\theta_{ij} - 1} \\
= & \int^{orb}_{\mc X_{(\mb g)}} e_1^*(\alpha) e_2^*(\beta) e_3^*(\gamma) \Theta_{(g_1,g_2)}\\
= & \int^{orb}_{\mc X_{(\mb g)}} e_1^*(\alpha) e_2^*(\beta) e_3^*(\gamma) e(E_{(\mb g)}).
\end{split}
\end{equation*}
Here $\Omega(X_{(\mb g)})$ is the Thom form 
of $X_{(\mb g)}$ in $X$, which represents the Poincar\'e dual of $X_{(\mb g )}$.
The theorem then follows from definition of Chen-Ruan orbifold cup product (\S \ref{ab:prod}).
\qed

We therefore constructed a deRham type model of $H^\ast_{CR}(\mc X)$. 
The advantages with this formulation is two-fold. Firstly, the product on Chen-Ruan orbifold cohomology can now be given directly. Secondly, as shown in the proof of the theorem, when computing the three point functions, the domain of integration are unified to be $\mc X$. This will make it easier for application.

\section{Symplectic reduction for torus action and wall crossing}
\label{wall}

As an application of our deRham model
 we consider symplectic reduction for torus action. 
Let $G = T^l$, $(M, \omega)$ be a $2N$-dimensional symplectic manifold with a 
Hamiltonian $G$ action. Let the moment map be $\mu : M \to \mathfrak{g}^*$, 
$p \in \mathfrak{g}^*$ lying in the image of $\mu$ be a regular value and 
$M(p) = \mu^{-1}(p)$. Then it's well known that 
$ X_p = M //_p G = M(p) / G$ is a symplectic orbifold 
of dimension $2n = 2(N-l)$. It is known that there is a chamber structure on 
$\mathfrak{g}^*$ such that 
$X_p$ and $X_q$ are diffeomorphic when $p$ and $q$ are in a same chamber $C$.
It would be interesting to investigate how the orbifold cohomology differs
when $p$ and $q$ are in different chambers. In this section, we will give
a wall crossing formula for the 3-point function. As one expects, the 
difference of 3-point functions on $X_p$ and $X_q$ is contributed by
fixed loci of the $G$ action on $M$. With the original formulation 
given in \S\ref{ab}, it is not easy to write a clean wall crossing formula 
due to the appearance of twisted sectors and obstruction forms. The new formulation
then has an advantage in dealing with these issues, at least at the level of
presentation.

\subsection{Orbi-structure of $X_p$}
\label{wall:str}

Let $\pi: M(p) \to X_p$ be the quotient map. Let $x \in X_p$ and $\tilde{x} \in \pi^{-1}(x)$, then a local orbifold chart $U$ near $x$ 
is given by a normal slice at $\tilde{x}$ of the orbit $G\circ\{\tilde{x}\}$ in $M(p)$, where $G_x$ is the isotropy group at $\tilde{x}$. Since $G$ is abelian, we are in the situation discussed in \S\ref{ab}.
%
The local group $G_x$ is a finite subgroup of $G$ and we make the following non-essential assumption to simplify notations. \footnote{Otherwise the labeling set $T_k$ below would have to take into account different components of the points fixed by subgroup $H$ since the same elements in local groups for different component are \emph{not} equivalent, which only leads to messier notations.}
\begin{assumption}
\label{singlecomponent}
For any finite subgroup $H$ of $G$, the points of $M$ which are fixed by $H$ have at most one component over any $p \in \mathfrak{g}^*$.
\end{assumption}

Under this assumption, the labeling set $T^k$ for $k$-multisectors of $X_p$ is subset of $G^k$. Let $\mathbf{g} = (g_1, \dots, g_k) \in G_x^k$, $\ggroup$ be the subgroup of $G_x$ generated by $\mathbf{g}$, $M^{\ggroup}$ be the fixed point set of $\ggroup$ in $M$ and $M(p)^{\ggroup} = M^{\ggroup} \cap M(p)$. Then we have
$$\(X_p\)_{(\mb g)} = M(p)^{\ggroup} / G.$$

\subsection{Equivariant set-up on $M$}
\label{wall:equi}

For simplicity, we will assume $G=S^1$.
Let $\mathcal{F}$ be the set of fixed points of $S^1$.
For $g\in G$, define $M^g$ to be the submanifold in $M$ fixed by
$g$. The interesting case is that $M^g-\mc F\not=\emptyset$.
From now on, we always assume that this is the case. 

The $G$ action gives a $G$-equivariant decomposition
$$TM|_{M^g} = \oplus_{j = 1}^m \tilde{E}_j \oplus TM^{g }.$$
This decomposition descends to the one in \eqref{ab:tandecom} with $TM^g$ further splits into $\mathbb R \oplus TM(p)^g
$ on $M(p)$. 
Let $[\tilde{l}_j]$ be an equivariant Thom class for $\tilde{E}_j$ supported in an equivariant neighbourhood of $0$-section of $\tilde{E}_j$. Let $\theta_j$ be the weights of $g$ action on fiber of $\tilde{E}_j$ then
\begin{defn}
\label{wall:fact}
\emph{Equivariant twist factor} for $M^g$ is the formal equivariant form:
$$[\tilde{t}(g)] = \prod_{j = 1}^m[ \tilde{l}_j]^{\theta_j}.$$
\end{defn}
As before, formally we have $\tilde{t}{(g)} \in H^{2 \iota_{(g)}}_G(M)$. We then make the following definitions parallel to those in \S\ref{dr:fact}
$$\tilde{i}_{(g)}: H^*_G(M^g) \to H^*_G(M) : \tilde{\alpha} \mapsto \tilde{i}_{(g)}(\tilde{\alpha}) = \tilde{\pi}^*(\tilde{\alpha}) \tilde{t}{(g)}, \text{ and }$$
$$H^*_{G,CR}(M) = \oplus_{(g) \in G} H^*_G(M^g)$$
with the degree shifting given by $2\iota_{(g)}$.

The Kirwan map for the usual (equivariant) cohomology is defined for regular value $p$ of the moment map $\mu$ as following:
$$\kappa_p : H^*_G(M) \xto{i_p^*} H^*_G(M(p)) \xto{\cong} H^*(X_p).$$
Kirwan surjectivity (\cite{K}) states that $\kappa_p$ is surjective when $M$ is compact. For some cases of non-compact $M$, e.g. $\bb{C}^n$ with linear actions, the Kirwan map is also surjective.
Suppose $\kappa_p$ is surjective for $M$ as well as $M^{g}$ for all $g \in G$ and define 
$$\kappa_p: H^*_{G,CR}(M) \to H^*_{CR}(X_p)$$
 by the direct sum on the factors, then the following is obvious:
\begin{prop}
\label{equiv:kirwan}
The Kirwan map $\kappa_p$ is surjective.
\qed
\end{prop}

\subsection{Wall crossing of Chen-Ruan orbifold cup product}
\label{wall:prod}
The set of regular values of $\mu$ consists of points outside a collection of hyperplanes in $\mathfrak{g}^*$. The codimension $1$ hyperplanes are called  \emph{walls} of the moment map. Let $W$ be such a wall and let $\xi_1 \in \mathfrak{g}$ be a primary vector such that 
$$W \subset \{v \in \mathfrak{g}^* | \langle \xi_1, v\rangle = 0\}.$$
Extend $\xi_1$ to a basis $\{\xi_1, \ldots, \xi_l\}$ of $\mathbb Z^l$-lattice of $\mathfrak{g}$ and fix the basis in the following. Let $H$ be the subgroup generated by $\xi_1$ and $H'$ be generated by $\{\xi_2, \ldots, \xi_l\}$ be its complement. Let $\{u_i\}$ be the dual basis of $\{\xi_i\}$. Suppose $p \in \mathrm{Image}(\mu)$ be a regular value and $a \in \mathbb R^+$ small such that $q= p+a u_1$ are in different chambers separated by $W$. Let $I = [p, q]$ denote the line segment between $p$ and $q$.

Let $X_q = M//_q G = M(q)/G$ and $\tilde{M} = \mu^{-1}(I)$. Let $\mathbf{g} = (g_1, g_2, g_3)$ such that $g_1g_2g_3 = 1$ and $\alpha_p \in H^*(X_{(g_1)})$, $\beta_p \in H^*(X_{(g_2)})$ and $\gamma_p \in H^*(X_{(g_3)})$ with $\tilde{\alpha}$, $\tilde{\beta}$ and $\tilde{\gamma} \in H^*_{G,CR}(M)$ be their equivariant lifting respectively. Let $\alpha_q = \kappa_q(\tilde{\alpha})$ and so on. We may arrange them into the following diagram:
\begin{equation*}
\text{
\xymatrix{
& H^*_{G, CR}(M) \ar[dl]_{\kappa_p} \ar[dr]^{\kappa_q} & \\
H^*_{CR}(X_p) & & H^*_{CR}(X_q)
}
}
\end{equation*}
Let $F = \cup_j F_j \subset \tilde{M}$ be the fixed point set of $H$ action. Then our main theorem is that 
the contribution to the difference between $\langle \alpha_p\cup \beta_p,\gamma_p\rangle$ and 
$\langle \alpha_q\cup \beta_q,\gamma_q\rangle$ 
is localized at $F$.  To simplify notations, we state theorem for 
the case $G=H=S^1$.
\begin{theorem}
\label{prod:thm}
Suppose $G=S^1$, $\alpha_p,\beta_p,\gamma_p$,
$\alpha_q,\beta_q,\gamma_q$ and $F_j$ are given as above.
Then
\begin{equation} 
\label{loc}
\langle \alpha_q\cup \beta_q,\gamma_q\rangle
-\langle \alpha_p\cup \beta_p,\gamma_p\rangle
=
\sum_j
\int_{F_j} \frac{\tilde{i}_{(g_1)}(\tilde{\alpha}) \tilde{i}_{(g_2)}(\tilde{\beta})
\tilde{i}_{(g_3)}(\tilde{\gamma})}{e_G(N_{F_j})},
\end{equation}
where $e_G(N_{F_j})$ is the equivariant euler class
of normal bundle $N_{F_j}$ of $F_j$ in $M$.
\end{theorem}
\n
{\bf Proof.}
Note that, by theorem \ref{isothm}, we have
$$
\langle \alpha_r\cup \beta_r,\gamma_r\rangle
=\int^{orb}_{X_r}
i_{(g_1)}(\alpha_r) i_{(g_2)} (\beta_r) i_{(g_3)} (\gamma_r)
$$
for $r=p$ or $q$. Then \eqref{loc} follows from the 
standard localization formula for the integration
$$
\int_{\tilde M}
\tilde{i}_{(g_1)}(\tilde{\alpha}) \tilde{i}_{(g_2)}(\tilde{\beta})
\tilde{i}_{(g_3)}(\tilde{\gamma}).
$$
\qed

\section{Examples}
\label{appl}
\subsection{Weighted projective spaces}
\label{wtproj}
Weighted projective space $\bb{P}(w_1, \ldots, w_n)$ of (complex) dimension $n-1$ can be described as the symplectic quotient of a linear $S^1$ action $\rho$ on $M =\bb{C}^n$, with weights $w_1, \ldots, w_n \in \bb{Z}^{+}$ on the eigenspaces. Let $W = (w_1, \ldots, w_n) \in \bb{Z}^n$ be the weight vector of the $S^1$ action. For simplicity, we assume that the greatest common divisor of $w_i$'s is $1$. Let the basis $\{v_i\}$ of $\bb{C}^n$ be given by the eigenvectors, then the moment map is given by 
$$\mu: \bb{C}^n \to \bb{R} : \mu(z_1, \ldots, z_n) = \frac{1}{2}\sum_{i} w_i |z_i|^2.$$
$0$ is the only singular value which is the wall and $\mu^{-1}(0) = \{0\} \in \bb{C}^n$ is the unique fixed point. Let $p < 0$ and $q > 0$ then we have $X_p = \emptyset$ and $X_q = \bb{P}(w_1, \ldots, w_n)$ with scaled symplectic (or K\"ahler) form. It follows that in \eqref{loc} there is only one term on either side and $e_G(N_F)$ here is simply $u^n\prod_i w_i$. Furthermore, the twisted sectors are copies of lower dimensional weighted projective subspaces with weights $(w_{i \in I})$ for some $I \subset \{1, \ldots, n\}$ and we denote them $\bb{P}_I(W)$. Thus we have
\begin{equation}
\label{wproj}\langle \alpha_1\cup \alpha_2,\alpha_3\rangle = \left.\left(\frac{\tilde{i}_{(g_1)}(\tilde{\alpha}_1) \tilde{i}_{(g_2)}(\tilde{\alpha}_2)
\tilde{i}_{(g_3)}(\tilde{\alpha}_3)}{u^n\prod_i w_i}\right)\right|_{z=0},
\end{equation}
for $\alpha_i \in X_{(g_i)} \cong \bb{P}_{I_i}(W)$ and $g_1g_2g_3 = 1$. The evaluation at $z=0$ implies that only the terms with no form part contribute in the various equivariant twisted forms. 

Let's apply the formula \eqref{wproj} to $X = \bb{P}(W)$ where $W = (1,2,2,3,3,3)$, which is studied in \cite{J}. Let $g = \omega$ the $3$-rd root of $1$, then the twisted sector $X_{(g)}$ of $X$ defined by $g$ is isomorphic to $\bb{P}(3,3,3)$, or equivalently, $\bb{P}^2$ with trivial $\bb{Z}_3$ action. It's straight forward to see that $\iota_{(g)} = \frac{1}{3} + \frac{2}{3} + \frac{2}{3} = \frac{5}{3}$. Let $\alpha_i \in H^*(X_{(g)})$ for $i = 1, 2, 3$, then in order for $\langle \alpha_1\cup \alpha_2,\alpha_3\rangle \neq 0$, we must have $\alpha_i \in H^0(X_{(g)})$. Without loss of generality, let $\alpha_i = 1_{(g)}$. Then applying \eqref{wproj} we have
\begin{equation*}
\begin{split}
&\langle 1_{(g)}\cup 1_{(g)},1_{(g)}\rangle = \left.\left(\frac{\(\tilde{t}_{(g)}\)^3}{u^5\prod_i w_i}\right)\right|_{z=0}, \\
=& \left.\left(\frac{\((u+\cdots)^{\frac{1}{3}} (2u+\cdots)^{\frac{2}{3}} (2u+\cdots)^{\frac{2}{3}}\)^3}{2^23^3u^5}\right)\right|_{z=0} \\
=& \frac{4}{27}
\end{split}
\end{equation*}
where $\cdots$ stands for terms evaluating to $0$ when $z=0$. This verifies the computation in \cite{J}.

\subsection{Mirror quintic orbifolds}\label{mirquint}
We consider the mirror quintic orbifold $Y$, which is defined as a generic member of the anti-canonical linear system in the following quotient of $\bb{P}^4$ by $(\bb{Z}_5)^3$:
$$[z_1:z_2:z_3:z_4:z_5] \sim [\xi^{a_1}z_1:\xi^{a_2}z_2:\xi^{a_3}z_3:\xi^{a_4}z_4:\xi^{a_5}z_5],$$
where $\sum a_i \equiv 0 \mod 5$ and $\xi = e^{\frac{2 \pi i}{5}}$. Let $\triangle^{\circ}$ be the polytope with vertices $v_0 = e_0= (-1, -1, -1, -1)$ and $v_i = e_0 + 5e_i$ for $i = 1, 2, 3, 4$ where $\{e_i\}_{i= 1}^4$ is the standard basis of $\bb{R}^4$. Then coning the faces of $\triangle^{\circ}$ gives the fan $\Sigma$ which defines $X = \bb{P}^4/(\bb{Z}_5)^3$. We can also obtain $Y$ as the quotient by the above $(\bb{Z})^3$ of a quintic of the following form:
$$\tilde{Y} = \{z_1^5 + z_2^5 + z_3^5 + z_4^5 + z_5^5 + \psi z_1 z_2z_3z_4z_5 = 0\}, \text{ where } \psi^5 \neq -5^5.$$
The computation for mirror quintic was first done in \cite{PP}.

The ordinary cup product on $Y$ is computed in \cite{PP} \S 6 and we refer to there for details. We also follow \cite{PP} \S 5 for the description of twisted sectors of $Y$. The twisted sectors of $Y$ are either points or curves. The main simplification in applying our method is to compute the contribution from twisted sectors which are curves. Let $Y_{(\mb{g})}$ be a triple twisted sector which is an orbifold curve, where $(\mb{g}) = (g_1, g_2, g_3)$. Such curve only occurs as intersection of $Y$ with some $2$-dimensional invariant variety of $X$. It follows then the isotropy group for generic point in $Y_{(\mb{g})}$ can only be $G \cong \bb{Z}_5$ and we have $g_i \in G$. Furthermore, under the evaluation maps to $Y$, $Y_{(g_i)}$ and $Y_{(\mb{g})}$ have the same images, which we'll denote as $Y_{(G)}$. 

Using the deRham model, we note that the formal maps 
$$i_{(\cdot)} : H^*(Y_{(\cdot)}) \to H^{* + \iota_{(\cdot)}}_{CR}(Y)$$ 
where $\cdot$ is one of $g_i$ or $\mb{g}$, all factor through a tubular neighbourhood of $Y_{(G)}$ in $Y$. Since $Y$ is orbifold Calabi-Yau, the degree shifting $\iota_{(\cdot)}$ is always non-negative integer. In particular, if $g_i \neq id \in G$, we have to have $\iota_{(g_i)} = 1$. Let $\alpha_i \in H^*(Y_{(g_i)})$ and we consider the Chen-Ruan cup product $\alpha_1 \cup \alpha_2$. It suffices to evaluate the non-zero pairing of the following form
\begin{equation}\label{triprod}
\langle \alpha_1 \cup \alpha_2, \alpha_3 \rangle = \int^{ord}_{Y}\wedge_{i = 1}^3i_{(g_i)}(\alpha_i) \neq 0.
\end{equation}
When $g_3 =id$, we see that the Chen-Ruan cup product reduces to (ordinary) Poincar\'e duality. When $g_i \neq id$ for $i =1, 2, 3$, by direct degree checking we find that $\alpha_i \in H^0(Y_{(g_i)})$ for all $i$ and the wedge product in \eqref{triprod} is a multiple of product of twist factors $t(g_i) = [l_1]^{\theta_{i1}} [l_2]^{\theta_{i2}}$, where $[l_j]$'s are the Thom classes of the line bundle factors of the normal bundle. Without loss of generality, let $\alpha_i = 1_{(g_i)}$. Since $g_i \neq id$ by assumption, we have $\theta_{ij} > 0$ for all $i,j$. Thus 
\begin{equation}
\label{reduce}\int^{ord}_{Y}\wedge_{i = 1}^3i_{(g_i)}(1_{(g_i)}) = \int^{orb}_{Y_{(G)}} c_{\cdot},
\end{equation}
where $c_{\cdot}$ stands for the Chern class corresponding to either $[l_1]$ or $[l_2]$. Let $X_2$ be the $2$ dimensional invariant subvariety of $X = \bb{P}^4/(\bb{Z}_5)^3$ such that $Y \cap X_2 = Y_{(G)}$. Then there are $2$ invariant subvarieties $X_{3,1}$ and $X_{3,2}$ of dimension $3$ which contains $X_2$. Let $Y_j = Y \cap X_{3,j}$ for $j = 1, 2$. Then $c_j$ above is simply the Chern class of the normal bundle of $Y_{(G)}$ in $Y_j$. To finish the computation, we note that the whole local picture can be lifted to $\tilde{Y} \subset \bb{P}^4$ where the Chern classes corresponding to $c_j$ obviously integrate to $5$. Then the quotient by $(\bb{Z}_5)^3$ gives the answer to the integration \eqref{reduce} as
$$\int^{orb}_{Y_{(G)}} c_{\cdot} = \frac{5}{125} = \frac{1}{25},$$
which verifies the computation in \cite{PP}.

\end{document}